\documentclass{article}
\usepackage{amssymb, amsmath,amsthm}

\renewcommand{\epsilon}{\varepsilon}
\renewcommand{\phi}{\varphi}

\begin{document}
\parindent0 in
\parskip 1 em
\title{On the Weil-Petersson Geometry of Calabi-Yau Moduli}
\author{Zhiqin Lu  and Eisuke Natsukawa \footnote{
The first
author is partially supported by  NSF Career award DMS-0347033 and the
Alfred P. Sloan Research Fellowship.}
}

\date{September 6, 2005}

\maketitle

Let $M$ be the moduli space of polarized Calabi-Yau manifolds. Let $\omega$ be the Weil-Petersson  metric on $M$. Let $f$ be an invariant polynomial of ${\rm Hom}\,(TM,TM)$.  Let $\Theta$ be the curvature tensor of $\omega$. Among the other results,  we proved that, for any non-negative integer $l$, 
\[
\int_M f(\Theta)\wedge\omega^l
\]
is finite and is a rational number. The  result generalized the result of ~\cite{SL-2}.

We also studied the incompleteness of the Weil-Petersson metric and found the necessary and sufficient condition of the incompleteness of the Weil-Petersson metric in several variables. The result generalized a result of Wang~\cite{Wang} in the one dimensional case.

\end{document}